# OPTIMALITY OF NEIGHBOR-BALANCED DESIGNS FOR TOTAL EFFECTS

By R. A. Bailey and P. Druilhet

*Queen Mary, University of London and CREST-ENSAI*

The purpose of this paper is to study optimality of circular neighbor-balanced block designs when neighbor effects are present in the model. In the literature many optimality results are established for direct effects and neighbor effects separately, but few for total effects, that is, the sum of direct effect of treatment and relevant neighbor effects. We show that circular neighbor-balanced designs are universally optimal for total effects among designs with no self neighbor. Then we give efficiency factors of these designs, and show some situations where a design with self neighbors is preferable to a neighbor-balanced design.

**1. Introduction.** In many experiments, especially in agriculture, the response on one plot may be affected by treatments on neighboring plots as well as by the treatment applied to that plot. Similarly, in cross-over designs the response on one subject in a given period may be affected by the residual effects of treatments applied to that subject in the previous periods. In both situations optimality results of designs with neighbor balance properties are available in the literature. See Shah and Sinha (1989) for a review of optimality results of neighbor-balanced cross-over designs with first-order residual effects, Hedayat and Afsarinejad (1978), Cheng and Wu (1980), Kunert (1984b), Kushner (1997) for cross-over designs, Kunert (1984a) and Druilhet (1999) for circular designs.

However, optimality results are almost always established for treatment and neighbor (or residual) effects separately. Usually, one of the aims of the experiment is to find a single treatment which can be recommended for use on larger spatial areas or over longer time periods than those used for individual treatments in the experiment: for example, a single variety of wheat to be grown in whole fields, a single drug to alleviate the symptoms of chronic asthma, or a single type of feed to be given to cows throughout a whole lactation. When the chosen treatment is in use, its only neighbor will be itself; thus the effect of most importance is the sum of the direct effect of the treatment and the neighbor effect(s) of the same treatment.

For cross-over trials, this is the sum of the direct effect and the residual effect of the same treatment. In the context of animal feeding trials, Patterson (1950, 1951) called this the *total effect* of the treatment; Kempton (1991) the *permanent effect*. Matthews (1988) also used the term *total effect* in cross-over trials, no matter what the area of application. In experiments on growing plants, the treatments might be varieties or species. In this case the sum of the direct effect of a treatment and its own neighbor effect is sometimes called the *species effect* [McGilchrist (1965)], the *variety effect* [Besag and Kempton (1986)], the *monoculture effect* [McGilchrist and Trenbath (1971) and Kempton (1982)] or the *pure stand effect* [Kempton (1985, 1991, 1997)]. However, the treatments may equally well be pesticides or quantities of fertilizer, in which case the French "l'effet plein champ" seems more appropriate. However, this terminology is clearly not suitable for drug trials. We propose that the sum of the direct effect of a treatment and any relevant neighbor effects of that treatment should be called the "total effect" irrespective of the area of application.

In some contexts it has been proposed that the effects of direct and neighbor effects are not additive when a treatment has itself as a neighbor [Speckel, Vincourt, Azaïs and Kobilinsky (1987) for competition between sunflowers]. In others [e.g., Kempton (1991)], there are correlations between neighboring responses. We do not assume such extra complications here.

Our designs are in blocks, usually incomplete. Ignoring neighbor effects, a block design is inefficient if any treatment occurs more than once in a block, so we usually insist that our designs are binary. Thus no treatment is ever a neighbor of itself. In particular, it is impossible to achieve orthogonality between direct effects and neighbor effects. Some designs optimal separately for the estimation of direct treatment and neighbor effects are known in which no treatment is preceded by itself. Examples of such designs are given by Hedayat and Afsarinejad (1978) for cross-over designs without preperiod, and by Magda (1980) for circular cross-over designs. Azaïs, Bailey and Monod (1993) give a catalog of circular neighbor-balanced designs with $t-1$ blocks of size $t$ or $t$ blocks of size $t-1$, where $t$ is the number of treatments. These latter designs have practical importance in field experiments. Indeed, neighbor balance implies that the number of replications is divisible by $t-1$ and experimenters rarely have the resources for $2(t-1)$ or more replications. In the present paper we aim to show that these designs are universally optimal [in Kiefer's (1975) sense] for total effects under models which incorporate one-sided or two-sided neighbor effects among the class of all designs with no treatment preceded by itself.

However, if there are a large number of blocks, then the gain from having self neighbors can outweigh the loss from nonbinarity. We use continuous block design theory to examine such designs and hence derive efficiency factors for binary neighbor-balanced designs.

**2. The designs and the models.** All the designs are assumed to be in linear blocks, with neighbor effects only in the direction of the blocks (say left-neighbor and right-neighbor effects). Because the effect of having no treatment differs from the neighbor effect of any treatment, we consider only designs with *border plots*, that is, designs with one plot added at each end of each block. The border plots receive treatments but are not used for measuring the response variables [see Langton (1990) for the importance of such designs]. The plots which are not on the borders are called *inner plots*. The *length* of a block is the number of its inner plots. We assume that all the designs are *circular*, that is, the treatment on a border plot is the same as the treatment on the inner plot at the other end of the block.

We denote by $\Omega_{(t,b,k)}$ the set of all circular designs with $t$ treatments and $b$ blocks of length $k$. We assume that $k \leq t$ and that there is an integer $\ell$ equal to $bk\{t(t-1)\}^{-1}$.

DEFINITION 1. A *circular binary block design* is a circular design which has each treatment at most once in the inner plots of each block.

DEFINITION 2. A *circular neighbor-balanced design* (CNBD) is a circular binary block design in $\Omega_{(t,b,k)}$ which is a balanced block design in the usual sense [Shah and Sinha (1989)] and such that for each ordered pair of distinct treatments there exist exactly $\ell$ inner plots which receive the first chosen treatment and which have the second one as right neighbor.

DEFINITION 3. A *circular design neighbor-balanced at distances* 1 *and* 2 (CNBD2) is a circular neighbor-balanced design such that for each ordered pair of distinct treatments, there exist exactly $\ell$ inner plots that have the first chosen treatment as left neighbor and the second one as right neighbor.

Here are two examples of circular designs neighbor-balanced at distance 1 and 2. The rows correspond to the blocks and the plots at each end of the blocks are the border plots.

$$\begin{bmatrix} 5 & | & 1 & 2 & 3 & 4 & 5 & | & 1 \\ 4 & | & 2 & 5 & 3 & 1 & 4 & | & 2 \\ 1 & | & 3 & 5 & 2 & 4 & 1 & | & 3 \\ 5 & | & 4 & 3 & 2 & 1 & 5 & | & 4 \end{bmatrix}, \quad \begin{bmatrix} 4 & | & 2 & 3 & 4 & | & 2 \\ 3 & | & 1 & 4 & 3 & | & 1 \\ 2 & | & 4 & 1 & 2 & | & 4 \\ 1 & | & 3 & 2 & 1 & | & 3 \end{bmatrix}.$$

For a design $d$, we denote by $d(i,j)$ the treatment assigned to plot $j$ of block $i$: in particular, $d(i,0)$ and $d(i,k+1)$ are the two treatments applied to the border plots of block $i$. The circularity condition implies that $d(i,0) = d(i,k)$ and $d(i,k+1) = d(i,1)$. We also denote by $Y_{i,j}$ the response on plot $j$ of block $i$. All the observations are assumed to be uncorrelated with common variance. We deal with two distinct models for the expectation,

$(\mathcal{M}1)$ $\quad \mathbb{E}(Y_{i,j}) = \beta_i + \tau_{d(i,j)} + \lambda_{d(i,j-1)}$ $\quad$ for $1 \leq i \leq b$ and $1 \leq j \leq k$,

$(\mathcal{M}2)$ $\quad \mathbb{E}(Y_{i,j}) = \beta_i + \tau_{d(i,j)} + \lambda_{d(i,j-1)} + \rho_{d(i,j+1)}$
$\quad$ for $1 \leq i \leq b$ and $1 \leq j \leq k$.

The unknown parameters have the following meanings: $\beta_i$ is the effect of block $i$, $\tau_{d(i,j)}$ is the effect of treatment $d(i,j)$, $\lambda_{d(i,j-1)}$ is the left-neighbor effect of treatment $d(i,j-1)$, $\rho_{d(i,j+1)}$ is the right-neighbor effect of treatment $d(i,j+1)$. In the standard vector notation, we have

$(\mathcal{M}1)$ $\qquad \mathbb{E}(Y) = B\beta + T_d\tau + L_d\lambda,$

$(\mathcal{M}2)$ $\qquad \mathbb{E}(Y) = B\beta + T_d\tau + L_d\lambda + R_d\rho,$

where $B, T_d, L_d$ and $R_d$ are the incidence matrices of block, treatment, left-neighbor and right-neighbor effects. Note that model $(\mathcal{M}1)$ corresponds to only one-sided neighbor effect. It is particularly adapted to temporal problems with carry-over effects. Model $(\mathcal{M}2)$ corresponds to specific additive influence from each left and right neighbor.

DEFINITION 4. The vector $\phi$ of total effects for models with one-sided neighbor effects is defined by $\phi = \tau + \lambda$.

DEFINITION 5. The vector $\psi$ of total effects for models with two-sided neighbor effects is defined by $\psi = \tau + \lambda + \rho$.

**3. Some technical tools.** We introduce some notation and results used throughout the next sections.

We denote by $\mathbb{1}_k$, $I_k$ and $J_k$, respectively, the vector of ones of length $k$, the $(k,k)$ identity matrix and the $(k,k)$ matrix of ones. For any matrix $A$, we denote by $A^+$ the Moore–Penrose inverse of $A$. The projection matrix onto the column span of matrix $A$ is denoted by $\mathrm{pr}_{(A)}$. Thus $\mathrm{pr}_{(A)} = A(A'A)^+A'$. We also define $\mathrm{pr}_{(A)}^\perp$ by $\mathrm{pr}_{(A)}^\perp = I - \mathrm{pr}_{(A)}$. Put $Q_k = \mathrm{pr}_{(\mathbb{1}_k)}^\perp = I_k - k^{-1}J_k$. For a square matrix $A$, we denote by $\mathrm{tr}(A)$ the trace of $A$. For two symmetric matrices $M$ and $N$, $M \leq N$ means that $N - M$ is a nonnegative definite matrix. A matrix is *completely symmetric* if it can be written as $aI + bJ$ for two scalars $a$ and $b$.

Consider the standard partitioned linear model:

$$(\mathcal{M}) \quad Y = A\alpha + B\beta + \varepsilon \quad \text{with } \mathbb{E}(\varepsilon) = 0 \text{ and } \text{Var}(\varepsilon) = \sigma^2 I.$$

The first lemma is classical [Kunert (1983)], the second, whose proof is in the Appendix, gives an upper bound of the information matrix: it generalizes a result by Pukelsheim [(1993), page 97] and gives a simple condition for equality.

LEMMA 1. *Under model ($\mathcal{M}$) the information matrix $C[\alpha]$ for the effect $\alpha$ is $C[\alpha] = A' \operatorname{pr}_{(B)}^{\perp} A$.*

LEMMA 2. *Assume that in model ($\mathcal{M}$) our interest is just for some linear combinations of $\alpha$, say $K'\alpha$. Then the corresponding information matrix $C[K'\alpha]$ satisfies $C[K'\alpha] \leq (K'K)^+ K' C[\alpha] K (K'K)^+$ with equality if and only if $C[\alpha]$ commutes with $\operatorname{pr}_{(K)}$.*

Throughout the paper we deal with universal optimality defined by Kiefer (1975). A universally optimal design has many good properties: see Shah and Sinha (1989) for further details.

PROPOSITION 3 [Kiefer (1975)]. *Assume that a design $d^*$ has its information matrix completely symmetric; then $d^*$ is universally optimal over a class $\mathcal{D}$ of designs if and only if $\operatorname{tr}(C_{d^*}) = \max_{d \in \mathcal{D}} \operatorname{tr}(C_d)$.*

**4. Optimality of circular neighbor-balanced designs for total effects.** In this section, we show that a CNBD (resp. a CNBD2) is universally optimal among all the designs with no treatment preceded by itself under the one-sided (resp. two-sided) neighbor effect model.

LEMMA 4. *Let $d^*$ be a circular neighbor-balanced design in $\Omega_{(t,b,k)}$ with $3 \leq k \leq t$; then the information matrix $C_{d^*}[\phi]$ for total effects under model ($\mathcal{M}1$) is $C_{d^*}[\phi] = b(k-2)(2(t-1))^{-1} Q_t$.*

LEMMA 5. *Let $d^*$ be a circular design neighbor-balanced at distances 1 and 2 in $\Omega_{(t,b,k)}$, with $4 \leq k \leq t$. Then the information matrix $C_{d^*}[\psi]$ for total effects under model ($\mathcal{M}2$) is $C_{d^*}[\psi] = b(k-3)(3(t-1))^{-1} Q_t$.*

THEOREM 6. *Under model ($\mathcal{M}1$) and for $3 \leq k \leq t$, a circular neighbor-balanced design in $\Omega_{(t,b,k)}$ is universally optimal for the total effects among all the designs with no treatment neighbor of itself.*

THEOREM 7. *Under model (M2) and for $4 \leq k \leq t$, a circular design neighbor-balanced at distance 1 and 2 in $\Omega_{(t,b,k)}$ is universally optimal for the total effects among all the designs with no treatment neighbor of itself at distance 1 or 2.*

PROOF OF LEMMA 4 AND THEOREM 6. Put $\alpha' = (\tau'|\lambda')$. Then $\phi = K'\alpha$, with $K = \mathbb{1}_2 \otimes I_t$.

By Lemma 1 we have $C_d[\alpha] = (T_d|L_d)' \operatorname{pr}^\perp_{(B)} (T_d|L_d)$ for any design. Here $K'K = 2I_t$, so, by Lemma 2 we have

$$C_d[\phi] \leq \tfrac{1}{4} K' C_d[\alpha] K$$

$$(1) \qquad = \tfrac{1}{4} \{ T'_d \operatorname{pr}^\perp_{(B)} T_d + T'_d \operatorname{pr}^\perp_{(B)} L_d + L'_d \operatorname{pr}^\perp_{(B)} T_d + L'_d \operatorname{pr}^\perp_{(B)} L_d \}$$

$$= \tfrac{1}{4} \{ 4 T'_d \operatorname{pr}^\perp_{(B)} T_d + T'_d L_d + L'_d T_d - 2 T'_d T_d \}.$$

The last equality comes from the fact that because of the circularity we have $T'_d T_d = L'_d L_d$ and $B' T_d = B' L_d$. So $\operatorname{pr}_{(B)} T_d = \operatorname{pr}_{(B)} L_d$ and so

$$T'_d \operatorname{pr}^\perp_{(B)} L_d = T'_d L_d - T'_d \operatorname{pr}_{(B)} L_d = T'_d L_d - T'_d T_d + T'_d \operatorname{pr}^\perp_{(B)} T_d.$$

Note that $(T'_d L_d)_{ii}$ is the number of times that treatment $i$ is a neighbor of itself at distance 1. Thus, for any design $d$ in $\Omega_{(t,b,k)}$ with no treatment neighbor of itself, $\operatorname{tr}(C_d[\phi])$ depends only on $\operatorname{tr}(T'_d \operatorname{pr}^\perp_{(B)} T_d)$. For a CNBD $d^*$ in $\Omega_{(t,b,k)}$, we have $T'_{d^*} L_{d^*} = L'_{d^*} T_{d^*} = \ell(J-I)$. Hence $C_{d^*}[\alpha]$ commutes with $\operatorname{pr}_{(K)} = 2^{-1} J_2 \otimes I_t$. Then Lemma 2 gives equality in (1). Because $T'_{d^*} \operatorname{pr}^\perp_{(B)} T_{d^*} = b(k-1)(t-1)^{-1} Q_t$, (1) establishes Lemma 4 for a CNBD. Moreover, a CNBD is a balanced block design, so it maximizes $\operatorname{tr}(T'_d \operatorname{pr}^\perp_{(B)} T_d)$ among all possible designs of the same size [see Shah and Sinha (1989)] and so, by Proposition 3 and Lemma 4, Theorem 6 is established. $\square$

PROOF OF LEMMA 5 AND THEOREM 7. Similarly, under model (M2), $\psi = K'\alpha$ with $\alpha' = (\tau'|\lambda'|\rho')$ and $K = \mathbb{1}_3 \otimes I_t$. Because of the circularity, we have $R'_d T_d = T'_d L_d$ and $T'_d R_d = L'_d T_d$. Then we have

$$(2) \quad C_d[\psi] \leq \tfrac{1}{9} \{ 9 T'_d \operatorname{pr}^\perp_{(B)} T_d + 2(T'_d L_d + L'_d T_d) + L'_d R_d + R'_d L_d - 6 T'_d T_d \}.$$

For any design with no treatment neighbor of itself at distance 1 or 2, $\operatorname{tr}(C_d[\psi])$ depends only on $\operatorname{tr}(T'_d \operatorname{pr}^\perp_{(B)} T_d)$. For a CNBD2 $d^*$ in $\Omega_{(t,b,k)}$, we have $R'_{d^*} L_{d^*} = L'_{d^*} R_{d^*} = \ell(J-I)$. Then Lemma 2 gives equality in (2). The remainder of the proof is identical to the previous one. $\square$

**5. Efficiency of circular neighbor-balanced designs.** In this section we examine efficiency factors of circular neighbor-balanced designs. First we use the method developed by Kushner (1997) and Kunert and Martin (2000) to construct optimal designs. Even if these designs usually have a large number of blocks, they are useful to derive efficiency factors.

5.1. *Continuous block designs and related upper bounds.* For each block $u$ of design $d$, we denote by $T_{du}$, $L_{du}$ and $R_{du}$ the incidence matrices corresponding to block $u$. Thus, $T_d = (T'_{d1}|\cdots|T'_{db})'$, $L_d = (L'_{d1}|\cdots|L'_{db})'$ and $R_d = (R'_{d1}|\cdots|R'_{db})'$. Inequalities (1) and (2) give, respectively,

$$(3) \qquad C_d[\phi] \leq \sum_{u=1}^{b} C_{du} \quad \text{and} \quad C_d[\psi] \leq \sum_{u=1}^{b} \tilde{C}_{du},$$

where $C_{du} = \frac{1}{4}\{4T'_{du}Q_kT_{du} + T'_{du}L_{du} + L'_{du}T_{du} - 2T'_{du}T_{du}\}$ and

$\tilde{C}_{du} = \frac{1}{9}\{9T'_{du}Q_kT_{du} + 2(T'_{du}L_{du} + L'_{du}T_{du}) + L'_{du}R_{du} + R'_{du}L_{du} - 6T'_{du}T_{du}\}$.

Because $\operatorname{tr}(C_{du})$ and $\operatorname{tr}(\tilde{C}_{du})$ are invariant under permutations of treatment labels, we may say that two sequences of treatments on a block are equivalent if one sequence can be obtained from the other one by relabelling the treatments. If we denote by $s$ the equivalence class of the sequence $l$ on the block $u$, we can define

$$c(s) = \operatorname{tr}(C_{du}) = \frac{1}{2}\left(k - \frac{2}{k}\sum_{i=1}^{t} n_i^2 + \sum_{i=1}^{t} m_i\right),$$

$$\tilde{c}(s) = \operatorname{tr}(\tilde{C}_{du}) = \frac{1}{9}\left(3k - \frac{9}{k}\sum_{i=1}^{t} n_i^2 + 4\sum_{i=1}^{t} m_i + 2\sum_{i=1}^{t} p_i\right),$$

where $n_i$ is the number of occurrences of treatment $i$ in the sequence $l$, $m_i$ is the number of times treatment $i$ is on the left-hand side of itself in sequence $l$ and $p_i$ is the number of plots having treatment $i$ on the left-hand side and the right-hand side.

PROPOSITION 8. *Consider a design $d$ with $b$ blocks of size $k$. If $s^*$ maximizes $c(s)$, then, under model $(\mathcal{M}1)$, $\operatorname{tr}(C_d[\phi]) \leq bc(s^*)$. If $s^*$ maximizes $\tilde{c}(s)$, then, under model $(\mathcal{M}2)$, $\operatorname{tr}(\tilde{C}_d[\psi]) \leq b\tilde{c}(s^*)$.*

PROOF. Denote by $S$ the number of equivalence classes and by $\pi_d(s)$ the proportion of the blocks containing a sequence in the equivalence class $s$. Then we have $\operatorname{tr}(C_d[\phi]) \leq b\sum_{s=1}^{S} \pi_d(s)c(s) \leq bc(s^*)$. □

5.2. *Optimal continuous block designs under model $(\mathcal{M}1)$.* Here we characterize optimal sequences under model $(\mathcal{M}1)$ and show how to construct universally optimal designs without restriction on the competing classes of designs. Then we calculate an efficiency factor for a CNBD.

NOTATION 1. We denote by $\lfloor x \rfloor$ the integer part of the real $x$ and by $f$ the real function $f(v) = -1 + k - v/2 - (2 - v/k)\lfloor k/v \rfloor + v/k\lfloor k/v \rfloor^2$.

TABLE 1

| $k$ | 3 | 4 | 4 | 4 | 5 | 6 | 7 | 8 | 9 | 10 | 11 | 12 | 12 | 12 | 13 | 14 | 15 | 16 | 16 |
|---|---|---|---|---|---|---|---|---|---|---|---|---|---|---|---|---|---|---|---|
| $v^*$ | 3 | 2 | 3 | 4 | 3 | 3 | 4 | 4 | 4 | 5 | 5 | 4 | 5 | 6 | 5 | 5 | 5 | 5 | 6 |
| $v_-$ | 3 | 2 | 2 | 4 | 1 | 3 | 1 | 4 | 3 | 5 | 4 | 4 | 3 | 6 | 2 | 1 | 5 | 4 | 2 |
| $v_+$ | 0 | 0 | 1 | 0 | 2 | 0 | 3 | 0 | 1 | 0 | 1 | 0 | 2 | 0 | 3 | 4 | 0 | 1 | 4 |
| $n_-$ | 1 | 2 | 1 | 1 | 1 | 2 | 1 | 2 | 2 | 2 | 2 | 3 | 2 | 2 | 2 | 2 | 3 | 3 | 2 |
| $n_+$ | 2 | 3 | 2 | 2 | 2 | 3 | 2 | 3 | 3 | 3 | 3 | 4 | 3 | 3 | 3 | 3 | 4 | 4 | 3 |

PROPOSITION 9. *Consider model* ($\mathcal{M}1$) *and designs with blocks of size $k$ ($k \geq 3$) and $t$ treatments. Then a sequence $l^*$ in the class $s^*$ that maximizes $c(s)$ is characterized by:*

1. *The number of different treatments present in the sequence $l^*$ is $v^*$, where $v^*$ maximizes $f(v)$ subject to $v \in \{2, \ldots, t\}$.*
2. *The number of occurrences of a treatment present in the sequence is either $n_- = \lfloor k/v^* \rfloor$ or $n_+ = n_- + 1$.*
3. *The number of treatments that occur $n_-$ times in $l^*$ is $v_-$, where $v_- = k - v^* \lfloor k/v^* \rfloor$.*
4. *The number of treatments that occur $n_+$ times in $l^*$ is $v_+$, where $v_+ = v^* - v_-$.*
5. *Every treatment present in $l^*$ has all its occurrences side by side.*

*Moreover, $c(s^*) = f(v^*) \leq k - \sqrt{2k}$ with equality if $(2k)^{1/2}$ is an integer.*

The proof is given in the Appendix.

REMARK. When $k = 2q(q+1)$ for some positive integer $q$, then $f(v)$ is maximum at exactly three points: $v^* = 2q$, $v^* = 2q + 1$, $v^* = 2q + 2$. Put $w = k/\lfloor (1 + (2k+1)^{1/2})/2 \rfloor$. If $w$ is an integer, then $f(v)$ is maximum at exactly $v = w$. If not, then $f(v)$ is maximum at one or both of the two integers either side of $w$. Moreover, it can be shown that, when $k$ is large, $v^* \sim \lfloor (2k)^{1/2} \rfloor$.

EXAMPLE. Table 1 gives the composition of $l^*$ depending on $k$ (when $t \geq v^*$).

For example, for $k = 5$ an optimal sequence contains $v^* = 3$ treatments. One treatment appears once and two treatments appear twice. So, for instance,

$$l^* = (a, b, b, c, c).$$

For $k = 14$ an optimal sequence contains five treatments: one treatment appears twice and four treatments appear three times. So an optimal sequence

TABLE 2

| $k$ | 3 | 4 | 5 | 6 | 7 | 8 | 9 | 10 | 11 | 12 | 13 | 14 | 15 |
|---|---|---|---|---|---|---|---|---|---|---|---|---|---|
| Eff($d^*$) | 1 | 1 | 0.88 | 0.80 | 0.80 | 0.75 | 0.75 | 0.73 | 0.72 | 0.71 | 0.70 | 0.69 | 0.68 |

is, for instance,

$$(a,a,b,b,b,c,c,c,d,d,d,e,e,e).$$

Note that other optimal sequences can be deduced by circular permutation or symmetry. For $k=4$ or $k=12$, there are three possibilities for $l^*$ as seen in the remark above. For $k=16$ there are two possibilities.

THEOREM 10. *Consider designs with $b$ blocks of size $k$, $k \geq 3$, and $t$ treatments. Denote by $s^*$ an optimal equivalence class of sequences. Then a design $d^*$ that has each sequence in $s^*$ equally often is universally optimal among all possible designs with the same size. Moreover, $\text{tr}(C_{d^*}) = bf(v^*) = bc(s^*)$ where $v^*$ is the number of treatments present in one sequence of $s^*$.*

PROOF. By construction, all the sequences in $s^*$ are obtained from one sequence in $s^*$ by relabelling the treatments. Thus, for design $d^*$, $C_{d^*}$ is completely symmetric. Moreover, $T'_{d^*}L_{d^*} = L'_{d^*}T_{d^*}$. So, by the proof of Lemma 4, the condition for equality in Lemma 2 holds and we have $\text{tr}(C_{d^*}) = bc(s^*)$. So by Proposition 8, $d^*$ maximizes the trace and so $d^*$ is universally optimal. □

Consider now the classical criteria $\Phi_p(C_d) = \text{tr}(C_d^{-p}/(t-1))^{1/p}$, with $\text{tr}(M^q) = \text{tr}(M^+)^{-q}$ if $q < 0$, and $\Phi_0 = \lim_{p \to 0} \Phi_p$. It is well known that $\Phi_0, \Phi_1, \Phi_\infty$ correspond, respectively, to D-, A-, E-optimality [Shah and Sinha (1989)]. Moreover, $\Phi_{-1}(C_d) = 1/\text{tr}(C_d)$. For a completely symmetric matrix $C$, $\Phi_p(C)$ does not depend on $p$. Thus we can derive efficiency factors for any $\Phi_p$ of a CNBD $d^*$ relative to a continuous block design $d^{**}$ constructed in Theorem 10 by considering Eff($d^*$) = $\text{tr}(C_{d^*})/\text{tr}(C_{d^{**}})$, as Table 2 shows.

When $k$ is large, the efficiency factor for a CNBD in $\Omega_{(t,b,k)}$ can be approximated by $(k-2)(2(k-\sqrt{2k}))^{-1}$, which tends to 0.5 when $k$ tends to $+\infty$. For $k=3$ and 4, the efficiency factor of a CNBD is 1, so we have the following result:

PROPOSITION 11. *For $k=3$ or 4 a CNBD is universally optimal for total effects among all possible designs with equal size.*

Table 3

| $k$ | 4 | 5 | 6 | 7 | 8 | 9 | 10 | 11 | 12 | 13 | 14 |
|---|---|---|---|---|---|---|---|---|---|---|---|
| Eff$(d^*)$ | 1 | 0.88 | 0.75 | 0.7 | 0.65 | 0.60 | 0.59 | 0.57 | 0.55 | 0.54 | 0.53 |

5.3. *Efficiency of CNBD2 under model* ($\mathcal{M}2$). In this section we derive optimal sequences under model ($\mathcal{M}2$). Unlike the previous section, we cannot construct optimal designs from an optimal sequence by considering all the treatment relabellings of the initial sequence, essentially because the condition for equality in Lemma 2 does not hold. So we just indicate the main result without giving the proof and we derive an upper bound for the efficiency factors of a CNBD2.

NOTATION 2. We denote by $\tilde{f}$ the function
$$\tilde{f}(v_1, v_2) = -1 + k - 2v_1/3 - 8v_2/9$$
$$- (2 - 2v_1/k - v_2/k)\lfloor (k - v_1)/v_2 \rfloor + v_2/k\lfloor (k - v_1)/v_2 \rfloor^2.$$

PROPOSITION 12. *A sequence $l^*$ in an optimal equivalence class $s^*$ is characterized by maximizing $\tilde{f}(v_1, v_2)$, for all possible values attainable where $v_1$ is the number of treatments appearing once in $l^*$, and $v_2$ is the number of treatments appearing at least twice in $l^*$.*

REMARK. It can be shown that if $l^*$ is an optimal sequence, then not only must $v_1$ and $v_2$ have values $v_1^*$ and $v_2^*$ which maximize $\tilde{f}$, but also all the plots receiving the same treatment are placed side-by-side. For $k = 4$ an optimal sequence contains four different treatments and thus a CNBD2 is universally optimal among all possible designs. For $k \geq 6$, an optimal sequence does not contain any treatment just once (i.e., $v_1^* = 0$). As in Section 5.2, we can derive the efficiency factor of a CNBD2 (Table 3).

## APPENDIX

PROOF OF LEMMA 2. Put $P = \mathrm{pr}_{(K)} = K(K'K)^+K'$ and $M = I - P$. Then $\alpha = P\alpha + M\alpha$ and $M^2 = M$. So $\mathbb{E}(Y) = AK(K'K)^+(K'\alpha) + AM^2\alpha + B\beta$. By Proposition 2.3 of Kunert (1983), with his $A_d$, $B_{1d}$ and $B_{2d}$ replaced, respectively, by $AK(K'K)^+$, $B$ and $AM$, we have $C[K'\alpha] \leq (K'K)^+K'A'\mathrm{pr}_{(B)}^\perp AK(K'K)^+ = (K'K)^+K'C[\alpha]K(K'K)^+$ with equality if and only if $(K'K)^+K'A'\mathrm{pr}_{(B)}^\perp AM = 0$, or equivalently $(K'K)^+K'C[\alpha]M = 0$. If $C[\alpha]$ commutes with $P$, then $C[\alpha]$ commutes with $M$, so $(K'K)^+ K'C[\alpha]M = (K'K)^+ K'MC[\alpha] = 0$ because $K'M = 0$. Conversely, if $(K'K)^+ K'C[\alpha]M = 0$, then $PC[\alpha]M = 0 = MC[\alpha]P$, so $PC[\alpha] = PC[\alpha](P+M) = (P+M)C[\alpha]P = C[\alpha]P$. $\square$

PROOF OF PROPOSITION 9. Let $l$ be a sequence in $s$ and denote by $v$ the number of treatments present in $l$. If $v = 1$, then $c(s) = 0$; thus the maximum must be sought on $\{2, \ldots, t\}$. If $v > 1$, then necessarily $m_i \leq n_i - 1$ for $i$ such that $n_i > 0$. Thus, $\sum_i m_i \leq k - v$ and then

$$c(s) \leq \frac{1}{2}\left(2k - v - \frac{2}{k}\sum_{i=1}^{t} n_i^2\right) = A \qquad \text{(say)}$$

with equality if $m_i = n_i - 1$ for all $i$ such that $n_i \geq 1$, that is, if all the plots containing the same treatment are side-by-side. Fix $v$. Then, because $\sum_{i=1}^{t} n_i = k$, $A$ is maximum if and only if $n_i = \lfloor k/v \rfloor$ or $n_i = \lfloor k/v \rfloor + 1$ for any treatment $i$ present in the sequence. So, necessarily, the number of treatments present $\lfloor k/v \rfloor$ times in the sequence is $v\lfloor k/v + 1 \rfloor - k$ and the number of treatments present $\lfloor k/v \rfloor + 1$ times is $k - v\lfloor k/v \rfloor$. Points 2–4 of the proposition are then established.

For such a sequence, $A = f(v)$ where the function $f$ is defined in Section 5.2. The function $f$ is continuous (in spite of the integer part). For each positive integer $p$, $f$ is linear on $[k/(p+1), k/p]$ with slope $(2p^2 + 2p - k)/(2k)$. The slope increases with $p$, so $f$ is concave. If there is a positive integer $q$ such that $k = 2q(q+1)$, then the slope is zero on $[2q, 2q+2]$ and any real number in this interval maximizes $f$. However, $v$ is an integer, so $f(v)$ is maximum at exactly three points: $v^* = 2q$, $v^* = 2q+1$ and $v^* = 2q+2$. Otherwise, the slope is never zero, so $f$ is maximized only at $k/\lfloor (1+(2k+1)^{1/2})/2 \rfloor = w$. If $w$ is an integer, then $f(v)$ is maximum when $v = w$. Otherwise, it is maximum at one or both of the integers on either side of $w$.

Now, if $k/v$ is an integer, then we have $f(v) = (2k - 2k/v - v)/2 = g(v)$ (say). Because $g$ is concave and $f$ is linear on intervals, then for all $v$, $f(v) \leq g(v)$. The maximum of $g$ is at $v = (2k)^{1/2}$. Thus, we have $f(v^*) \leq g(\sqrt{2k}) = k - \sqrt{2k}$. $\square$

## REFERENCES

AZAÏS, J.-M., BAILEY, R. A. and MONOD, H. (1993). A catalogue of efficient neighbour-designs with border plots. *Biometrics* **49** 1252–1261.

BESAG, J. and KEMPTON, R. A. (1986). Statistical analysis of field experiments using neighbouring plots. *Biometrics* **42** 231–251. MR876840

CHENG, C.-S. and WU, C. F. (1980). Balanced repeated measurements designs. *Ann. Statist.* **8** 1272–1283. [Correction (1983) **11** 349.] MR594644

DRUILHET, P. (1999). Optimality of neighbour balanced designs. *J. Statist. Plann. Inference* **81** 141–152. MR1718460

HEDAYAT, A. and AFSARINEJAD, K. (1978). Repeated measurements designs. II. *Ann. Statist.* **6** 619–628. [MR488527](MR488527)

KEMPTON, R. A. (1982). Adjustment for competition between varieties in plant breeding trials. *J. Agricultural Sci.* **98** 599–611.

KEMPTON, R. A. (1985). Spatial methods in field experiments. *Biometric Bulletin* **2** (3) 4–5.

KEMPTON, R. A. (1991). Interference in agricultural experiments. In *Proc. Second Meeting of the Biometric Society East/Central/Southern African Network*. Harare, Zimbabwe.

KEMPTON, R. A. (1997). Interference between plots. In *Statistical Methods for Plant Variety Evaluation* (R. A. Kempton and P. N. Fox, eds.) 101–116. Chapman and Hall, London.

KIEFER, J. (1975). Construction and optimality of generalized Youden designs. In *A Survey of Statistical Design and Linear Models* (J. N. Srivastava, ed.) 333–353. North-Holland, Amsterdam. [MR395079](MR395079)

KUNERT, J. (1983). Optimal design and refinement of the linear model with applications to repeated measurements designs. *Ann. Statist.* **11** 247–257. [MR684882](MR684882)

KUNERT, J. (1984a). Designs balanced for circular residual effects. *Comm. Statist. A—Theory Methods* **13** 2665–2671. [MR759242](MR759242)

KUNERT, J. (1984b). Optimality of balanced uniform repeated measurements designs. *Ann. Statist.* **12** 1006–1017. [MR751288](MR751288)

KUNERT, J. and MARTIN, R. J. (2000). On the determination of optimal designs for an interference model. *Ann. Statist.* **28** 1728–1742. [MR1835039](MR1835039)

KUSHNER, H. B. (1997). Optimal repeated measurements designs: The linear optimality equations. *Ann. Statist.* **25** 2328–2344. [MR1604457](MR1604457)

LANGTON, S. (1990). Avoiding edge effects in agroforestry experiments: The use of neighbour-balanced designs and guard areas. *Agroforestry Systems* **12** 173–185.

MAGDA, C. G. (1980). Circular balanced repeated measurements designs. *Comm. Statist. A—Theory Methods* **9** 1901–1918. [MR603164](MR603164)

MATTHEWS, J. N. S. (1988). Recent developments in crossover designs. *Internat. Statist. Rev.* **56** 117–127. [MR963525](MR963525)

MCGILCHRIST, C. A. (1965). Analysis of competition experiments. *Biometrics* **21** 975–985.

MCGILCHRIST, C. A. and TRENBATH, B. R. (1971). A revised analysis of competition experiments. *Biometrics* **27** 659–671.

PATTERSON, H. D. (1950). The analysis of change-over trials. *J. Agricultural Sci.* **40** 375–380. [MR48774](MR48774)

PATTERSON, H. D. (1951). Change-over trials (with discussion). *J. Roy. Statist. Soc. Ser. B* **13** 256–271. [MR48774](MR48774)

PUKELSHEIM, F. (1993). *Optimal Design of Experiments*. Wiley, New York. [MR1211416](MR1211416)

SHAH, K. R. and SINHA, B. K. (1989). *Theory of Optimal Designs. Lecture Notes in Statist.* **54**. Springer, New York. [MR1016151](MR1016151)

SPECKEL, D., VINCOURT, P., AZAÏS, J.-M. and KOBILINSKY, A. (1987). Etude de la compétition interparcellaire chez le tournesol. *Biom. Praxim.* **27** 21–43.

SCHOOL OF MATHEMATICAL SCIENCES  
QUEEN MARY, UNIVERSITY OF LONDON  
MILE END ROAD  
LONDON E1 4NS  
UNITED KINGDOM  
E-MAIL: r.a.bailey@qmul.ac.uk

CREST-ENSAI  
ECOLE NATIONALE DE LA STATISTIQUE  
ET DE L'ANALYSE DE L'INFORMATION  
35 170 BRUZ  
FRANCE  
E-MAIL: druilhet@ensai.fr